\documentclass[reqno,12pt]{amsart}
\usepackage{amsmath, amssymb, amsthm,amsfonts} 
\usepackage[english]{babel}
\usepackage{bbm}
\usepackage{graphicx}
\usepackage{url}
\usepackage{epstopdf}
\usepackage[a4paper,bindingoffset=0.5cm,left=2cm,right=2cm,top=2.5cm,bottom=2cm,footskip=.8cm]{geometry}

\usepackage{tikz}
\usetikzlibrary{arrows, automata,positioning,calc,shapes,decorations.pathreplacing,decorations.markings,shapes.misc,petri,topaths}
\usepackage{tkz-berge}
  \usepackage{pgfplots}
  \pgfplotsset{compat=newest}
  \usetikzlibrary{plotmarks}
  \usepackage{grffile}
\newlength\figureheight
  \newlength\figurewidth
\pgfplotsset{%
    tick label style={font=\scriptsize},
    label style={font=\footnotesize},
    legend style={font=\footnotesize},
         every axis plot/.append style={very thick}
}

\usepackage{rotating}
\usepackage{amsbsy,enumerate}
\usepackage{graphicx}
\usepackage{ccaption}
\usepackage{comment}
\usepackage{mathrsfs}

\newcommand{\s}{^\star}

\newcommand{\bs}{\boldsymbol}

\newcommand{\vb}{\vspace{5.2mm}}
\renewcommand{\hat}{\widehat}

\newcommand{\vt}{\vartheta}

\newcommand{\DQ}{\mathscr Q}

\allowdisplaybreaks

\newtheorem{lemma}{Lemma}

\newtheorem{theorem}{Theorem}
\newtheorem{remark}{Remark}

\makeatletter
\renewcommand{\fnum@figure}[1]{\textbf{\figurename~\thefigure}. }
\renewcommand{\fnum@table}[1]{\textbf{\tablename~\thetable}. }
\makeatother

\begin{document}

\title[The correlation function of a queue with L\'evy and MAP input]{The correlation function of a queue \\with L\'evy and Markov additive input}

\author{Wouter Berkelmans, Agata Cichocka, Michel Mandjes}

\begin{abstract}
Let $(Q_t)_{t\in{\mathbb R}}$ be a stationary workload process, and $r(t)$ the correlation coefficient of $Q_0$ and $Q_t$. 
In a series of previous papers (i) the transform of $r(\cdot)$ has been derived for the case that the driving process is spectrally-positive ({\sc sp}) or spectrally-negative ({\sc sn}) L\'evy, (ii) it has been shown that for {\sc sp}-L\'evy and {\sc sn}-L\'evy input  $r(\cdot)$ is positive, decreasing, and convex, (iii) in case the driving L\'evy process is light-tailed (a condition that is automatically fulfilled in the {\sc sn} case), the decay of the decay rate agrees with that of the tail of the busy period distribution.  

\noindent
In the present paper we first prove the conjecture that property (ii) carries over to spectrally two-sided L\'evy processes; we do so for the case the L\'evy process is reflected at 0, and the case it is reflected at $0$ and $K>0$. Then we focus on queues fed by Markov additive processes ({\sc map}s). We start by the establishing the counterpart of (i) for {\sc sp}- and {\sc sn}-{\sc map}s. Then we refute property (ii) for {\sc map}s: we construct examples in which the correlation coefficient can be (locally) negative, decreasing, and concave. Finally, in relation to (iii), we point out how to identify the decay rate of $r(\cdot)$ in the light-tailed {\sc map} case, thus showing that the tail behavior of $r(\cdot)$ does not necessarily match that of the busy-period tail; singularities related to the transition rate matrix of the background Markov chain turn out to play a crucial role here. 

\vb


\vb

\noindent
{\sc Affiliations.} 
Wouter Berkelmans and Agata Cichocka are with 
Department of Mathematics,
Vrije Universiteit Amsterdam,
De Boelelaan 1081a,
1081 HV Amsterdam,
the Netherlands.
Michel Mandjes is with Korteweg-de Vries Institute for Mathematics, University of Amsterdam, Science Park 904, 1098 XH Amsterdam, the Netherlands. He is also with E{\sc urandom}, Eindhoven University of Technology, Eindhoven, the Netherlands, and Amsterdam Business School, Faculty of Economics and Business, University of Amsterdam, Amsterdam, the Netherlands. His research is partly funded by NWO Gravitation project N{\sc etworks}, grant number 024.002.003.

\vb

\noindent
{\sc Acknowledgments.} 
The third author wishes to thank Lars N{\o}rvang Andersen (Aarhus) and Fran\c{c}ois Baccelli (Austin) for helpful interactions regarding the spectrally two-sided L\'evy case, Jevgenijs Ivanovs (Aarhus) and Zbigniew Palmowski (Wroc{\l}aw) for providing useful references, and Peter Glynn (Stanford) for discussions regarding the tail behavior examples in Section \ref{maie}.

\end{abstract}

\maketitle

\newpage

\section{Introduction} 
Consider a queueing resource with workload process $(Q_t)_{t\in{\mathbb R}}$. Where much effort has been spent on characterizing transient and stationary distributional properties pertaining to this workload process, this paper focuses on the {\it workload correlation}
\[r(t) := \frac{{\mathbb C}{\rm ov}(Q_0,Q_t)}{\sqrt{{\mathbb V}{\rm ar}\,Q_0\cdot {\mathbb V}{\rm ar}\,Q_t}},\]
assuming the workload is in stationarity at time $0$. Insight into $r(\cdot)$ is useful in various ways, as it tells us from what time horizon $\tau$ on the random variables $Q_0$ and $Q_\tau$ can safely be assumed independent (in that $r(\tau)$ drops below some predefined level $\varepsilon$). 

In a series of articles, the function $r(\cdot)$ has been analyzed for the important class of {\it L\'evy-driven queues}. These can be seen as queues fed by L\'evy input \cite{DM}, and are also known as {\it reflected L\'evy processes} \cite{KY}.  Important subclasses of L\'evy processes are {\it spectrally positive} (in this paper abbreviated to {\sc sp}) L\'evy processes, with jumps only in the upward direction, and {\it spectrally negative}~({\sc sn}) L\'evy processes, with jumps only in the downward direction. The results on the workload correlation function $r(\cdot)$ that were derived so far can be divided into three types: (i) expressions for the transform $\hat r(\cdot)$ of $r(\cdot)$, (ii) structural properties concerning the shape of $r(\cdot)$, and (iii) characterizations of $r(t)$ for $t$ large. 

We proceed by giving a brief account of the literature. Ott \cite{OTT} considered a subclass of the {\sc sp} L\'evy input case, viz.\ queues with compound Poisson input (i.e., queues of the M/G/1 type). He succeeded in finding $\hat r(\cdot)$ in closed-form. In addition, he used this transform to prove, relying on the machinery of completely-monotone functions, that $r(\cdot)$ is non-negative, non-increasing, and convex.  
These results were extended by Es-Saghouani and Mandjes \cite{EM} to the full class of {\sc sp} L\'evy inputs (thus also covering L\'evy processes with a Brownian component, as well as components with so-called `small jumps' in the upward direction). In addition, they identified the tail behavior of $r(\cdot)$; in particular, when the {\sc sp} L\'evy input is light-tailed, it was shown that its (exponential) decay rate matches that of the tail of the busy-period distribution. Glynn and Mandjes~\cite{GM} proved the counterparts of the results of \cite{EM} for the case of {\sc sn} L\'evy input, and in addition devised a coupling-based simulation technique to efficiently estimate $r(t)$. The findings of \cite{EM,GM} show that the derivation of the Laplace transform $\hat r(\cdot)$ requires the availability of an explicit expression for the Laplace transform of the stationary workload, as well as for the Laplace transform of the transient workload (i.e., after an exponentially distributed interval $T$, conditional on the workload at time $0$).

Extensions in various directions could be thought of. Two of them are:
\begin{itemize}
\item[$\circ$] In the first place, one could consider queues with {\it general} L\'evy input (queues fed by a {\it spectrally two-sided} L\'evy process, that is). A long outstanding conjecture is that 
the fact that $r(\cdot)$ 
is non-negative, non-increasing, and convex, which was proven for {\sc sp} and {\sc sn} L\'evy input, carries 
over to the spectrally two-sided L\'evy setting. There is a  similar conjecture for the correlation of the waiting times in the GI/G/1 queue. 

\noindent
It is clear that a proof for the spectrally two-sided case requires a different approach than the one used in the one-sided cases. 
The major complication is that, unlike in the {\sc sp} and {\sc sn} cases, no explicit results  for the Laplace transforms of the stationary and transient workload are available --- results are typically in terms of `Wiener-Hopf expressions'; see e.g.\ \cite[Thm.\ 4.3]{DM}. As a result, identifying the transform $\hat r(\cdot)$ seems to be out of reach. As in the {\sc sp} and {\sc sn} cases this Laplace transform was crucial in establishing that $r(\cdot)$ 
is non-negative, non-increasing, and convex (using the concept of completely-monotone functions),
one may  wonder whether these structural properties can be proven without knowing $\hat r(\cdot).$ 

\noindent
In addition, one may wonder whether the property that $r(\cdot)$ 
is non-negative, non-increasing, and convex remains valid if the L\'evy process is {\it doubly reflected}, i.e., reflected both from below at level $0$ and from above at level $K>0$ (corresponding with the workload process of a queue with finite buffer capacity $K$). 

\vspace{1mm}

\item[$\circ$] Secondly, one could consider queues fed by spectrally one-sided {\it Markov additive processes} ({\sc map}s). These processes, dating back to \cite{CI,NE}, can be seen as the Markov-modulated of L\'evy processes. Informally, a {\sc map} can be seen as a L\'evy process, but with its local behavior depending on the state of an external finite-state Markov process, usually referred to as the {\it background process}. Only for very specific cases, partial results on the workload correlation function are known; see e.g.\ the findings in \cite{KM}. 
\end{itemize}

We now describe the contributions of our paper. In the first place, we prove the conjecture that in the case of two-sided L\'evy input the workload correlation function is non-negative, non-increasing, and convex. We do so without relying on explicit expressions for the transform $\hat r(\cdot)$; instead we provide a clean and insightful proof that uses first principles only. The argument carries over to the correlation of waiting times in the GI/G/1 queue. We also show that in case we add reflection at $K>0$, the property that $r(\cdot)$ 
is non-negative, non-increasing, and convex remains true. 

We then consider {\sc map}-driven queues, focusing on the cases with {\sc sp}-{\sc map} and {\sc sn}-{\sc map} input. For both classes we succeed in deriving an explicit expression for the Laplace transform $\hat r(\cdot)$, making intensive use of earlier fluctuation-theoretic results for spectrally one-sided {\sc map}s. One may wonder whether also for {\sc map} input the  workload correlation function is non-negative, non-increasing, and convex; by a series of insightful examples, we show that this is not the case. We also study the tail behavior of $r(\cdot)$, with the main conclusion that in the light-tailed setting the decay rate does not match with the busy-period decay rate (unlike in the L\'evy case). More specifically, as it turns out, the rate of decay of $r(\cdot)$ follows by comparing the busy-period decay rate with singularities  related to the transition rate matrix of the background process.

\section{L\'evy input: structural properties of correlation function}
In this section we consider the workload process pertaining to a queue with L\'evy input. In the first subsection we consider the case of
reflection from below at level 0, whereas the second studies its finite-buffer counterpart in which reflection from above at level $K>0$ is added.

\subsection{Single-sided reflection} Let $(X_t)_{t\in{\mathbb R}}$ be a single-dimensional L\'evy process \cite{KY} such that ${\mathbb E}\,X_1<0.$ Then we define $(Q_t)_{t\in{\mathbb R}}$ as the associated workload process (or: the process $(X_t)_{t\in{\mathbb R}}$ reflected at 0) by
\begin{equation}
\label{UR}Q_t := X_t +\max\left\{Q_0,-\inf_{0\leqslant u\leqslant t}X_u\right\}.\end{equation}
We assume the workload process is in stationarity, meaning that (under our stability condition ${\mathbb E}\,X_1<0$) we can alternatively write
\begin{equation}
\label{CT}Q_t = \sup_{s\leqslant t} (X_t-X_s);\end{equation}
due to the stationarity, we in particular have that, for any $t$, $Q_t$ is distributed as $Q_0.$
For future reference, we also introduce for $s\leqslant t$,
\[X_{s,t}:= X_t-X_s,\:\:\:\:\: L_{s,t} := -\inf_{s\leqslant u\leqslant t}X_{s,t}.\]
In addition, we let ${\mathscr G}_s$ be $\sigma\{(X_{s,t})_{t\geqslant s}\}$, i.e., the sigma-algebra generated by the increments of $(X_t)_{t\in{\mathbb R}}$ after time $s$ (relative to $X_s$, that is). In the sequel, we also work with $\bar Q$, which has the stationary workload distribution (and is therefore distributed as $Q_0$), but is independent of the process $(X_t)_{t\in{\mathbb R}}$ (and hence also of the process $(Q_t)_{t\in{\mathbb R}}$). 

For $s\leqslant t$, 
\begin{equation}\label{R1}Q_t =X_{s,t} +\max\{Q_s,L_{s,t}\} = X_{s,t}+Q_s +\max\{0,L_{s,t}-Q_s\};\end{equation}
observe that (i) $Q_t$ is non-decreasing in $Q_s$, and (ii) $L_{s,t}$ is independent of $Q_s$. This relation implies that, for $s\leqslant t\leqslant u$,
\begin{align}\nonumber Q_u-Q_t = &\:X_{t,u} +\max\{0,L_{t,u}-Q_t\}\\
=&\:X_{t,u} +\max\{0, L_{t,u}-X_{s,t}-\max\{Q_s,L_{s,t}\}\};\label{R2}\end{align}
observe that (i) $Q_u-Q_t$ is non-increasing in $Q_s$, and (ii) $L_{s,t}$, $X_{s,t}$, $L_{t,u}$, and $X_{t,u}$ are independent of $Q_s$. 

We throughout assume that ${\mathbb V}{\rm ar}\,Q_0<\infty$, implying that ${\mathbb C}{\rm ov}(Q_0,Q_t)<\infty$. The object of study is
\[r(t) := \frac{{\mathbb C}{\rm ov}(Q_0,Q_t)}{\sqrt{{\mathbb V}{\rm ar}\,Q_0\cdot {\mathbb V}{\rm ar}\,Q_t}}= 
\frac{{\mathbb E}\,Q_0Q_t - ({\mathbb E}\,Q_0)^2}{{\mathbb V}{\rm ar}\,Q_0}.\]
Our main objective is to prove, using first principles, that $r(\cdot)$ is non-negative, non-increasing, and convex (which was shown for {\sc sp}-L\'evy and {\sc sn}-L\'evy input in \cite{EM,GM}).

The following lemma will play a pivotal role in our analysis; its proof is provided in the appendix.
\begin{lemma} \label{LEM} Let $A$ and $B$ be identically distributed non-negative random variables. Let $f(\cdot)$ be non-decreasing, and $g(\cdot)$ be non-increasing. Then,
\[{\mathbb E}(A\,f(A)) \geqslant {\mathbb E}(B\,f(A)),\:\:\:{\mathbb E}(A\,g(A)) \leqslant {\mathbb E}(B\,g(A)).\]
\end{lemma}

The main result of this section is the following. 

\begin{theorem}\label{MT} In the system with the L\'evy process $(X_t)_{t\in{\mathbb R}}$ being reflected from below at 0, 
$r(\cdot)$ is non-negative, non-increasing, and convex.
\end{theorem}

{\it Proof}: (A)~We first show that ${\mathbb C}{\rm ov}(Q_0,Q_t)$ is non-negative. To this end, observe that
\begin{align*}
{\mathbb C}{\rm ov}(Q_0,Q_t) &\:= {\mathbb E}((Q_0 - \bar Q) Q_t) = {\mathbb E}\big({\mathbb E}((Q_0 - \bar Q) Q_t)\,|\,{\mathscr G}_0\big)\\
&\:= {\mathbb E}\big({\mathbb E}((Q_0 - \bar Q) (X_{0,t} +\max\{Q_0,L_{0,t}\} ))\,|\,{\mathscr G}_0\big),
\end{align*}
where the last equality is due to (\ref{R1}). Now apply Lemma \ref{LEM} with $A:=Q_0$, $B:=\bar Q$, and
\[f(a):= X_{0,t} +\max\{a,L_{0,t}\}.\]
To verify that the conditions  of the lemma are met, first note that $f(\cdot)$ is non-decreasing. In the second place, $Q_0$ is independent of ${\mathscr G}_0$ (i.e., the L\'evy process after time 0), and has therefore the same distribution as $\bar Q$. The lemma thus yields
\[{\mathbb C}{\rm ov}(Q_0,Q_t) = {\mathbb E}\big({\mathbb E}((Q_0 - \bar Q) (X_{0,t} +\max\{Q_0,L_{0,t}\} ))\,|\,{\mathscr G}_0\big)\geqslant 0,\]
as desired.

(B)~We now prove that ${\mathbb C}{\rm ov}(Q_0,Q_t)$ is non-increasing, for which it is sufficient that ${\mathbb E}(Q_0Q_t)$ is non-increasing. Again we use Lemma 1, but now with $A:=Q_s$, $B:=Q_0$ (with $0\leqslant s\leqslant t$), and
\[f(a):= X_{s,t}+\max\{a,L_{s,t}\}.\]
Evidently, $f(\cdot)$ is non-decreasing and
conditional on ${\mathscr G}_s$ the random variables $Q_0$ and $Q_s$ have the same (conditional) distributions. The lemma thus yields
\[{\mathbb E}\big({\mathbb E}(Q_0(X_{s,t}+\max\{Q_s,L_{s,t}\}))|\,{\mathscr G}_s\big) \leqslant 
{\mathbb E}\big({\mathbb E}(Q_s(X_{s,t}+\max\{Q_s,L_{s,t}\}))|\,{\mathscr G}_s\big).\]
We obtain
\begin{align*}
{\mathbb E}(Q_0Q_t) &\:= {\mathbb E}(Q_0(X_{s,t}+\max\{Q_s,L_{s,t}\}))={\mathbb E}\big({\mathbb E}(Q_0(X_{s,t}+\max\{Q_s,L_{s,t}\}))|\,{\mathscr G}_s\big)\\
&\:\leqslant {\mathbb E}\big({\mathbb E}(Q_s(X_{s,t}+\max\{Q_s,L_{s,t}\}))|\,{\mathscr G}_s\big)= {\mathbb E}(Q_sQ_t) = {\mathbb E}(Q_0Q_{t-s}),
\end{align*}
with the last equality being valid due to stationarity. 

(C)~In the third part of the proof, we show that ${\mathbb C}{\rm ov}(Q_0,Q_t)$ is convex. We do so by proving the sufficient property, for $s\leqslant t\leqslant u$,
\begin{equation}\label{goal}
{\mathbb E}(Q_0(Q_u-Q_t)\geqslant {\mathbb E}(Q_s(Q_u-Q_t)) = {\mathbb E}(Q_0(Q_{u-s}-Q_{t-s})),
\end{equation}
where the last equality is due to the stationarity; combined with the monotonicity of ${\mathbb E}(Q_0 Q_t)$ this implies convexity. In this case we define
\[f(a) := X_{t,u}+\max\{0,L_{t,u}-X_{s,t}-\max\{a,L_{s,t}\}\}.\] 
Again we wish to apply Lemma \ref{LEM} with $A:=Q_s$ and $B:=Q_0$. Now $f(\cdot)$ is non-increasing. In addition, conditional on ${\mathscr G}_s$, we again have that $Q_0$ and $Q_s$ are identically distributed. It thus follows that, using (\ref{R2}),
\begin{align*}
{\mathbb E}(Q_0(Q_u-Q_t)) &\:= {\mathbb E}(Q_0(X_{t,u} +\max\{0, L_{t,u}-X_{s,t}-\max\{Q_s,L_{s,t}\}\}))\\
&\:={\mathbb E}\big( {\mathbb E}(Q_0(X_{t,u} +\max\{0, L_{t,u}-X_{s,t}-\max\{Q_s,L_{s,t}\}\}))|\,{\mathscr G}_s\big)\\
&\:\geqslant {\mathbb E}\big({\mathbb E}(Q_s(X_{t,u} +\max\{0, L_{t,u}-X_{s,t}-\max\{Q_s,L_{s,t}\}\}))|\,{\mathscr G}_s\big)\\
&\:={\mathbb E}(Q_s(Q_u-Q_t)).
\end{align*}
This completes the proof. \hfill$\Box$

\begin{remark} {\em Waiting times in the GI/G/1 queue are known to satisfy the {\it Lindley recursion}. Indeed, with $W_n$ the waiting time of the $n$-th customer, $A_n$ the time between the arrivals of the $n$-th and $(n+1)$-st arrivals, $S_n$ the service time of the $n$-th customer, and
$U_n:=  S_n -A_n$,
\[W_{n+1} = \max\{W_n+U_n,0\}.\]
Iterating this relation, we obtain that, under the stability condition ${\mathbb E}(S_n)<{\mathbb E}(A_n)$,
\[W_{n+1} =\max\left\{\max_{i=-\infty,\ldots,n}\sum_{j=i}^n U_j,0\right\},\]
which is the discrete-time counterpart of (\ref{CT}). Mimicking all steps that we used in the proof of Thm.\ \ref{MT}, it follows that the correlation function of $(W_n)_{n\in{\mathbb Z}}$ is non-negative, non-increasing and convex as well. $\hfill\diamond$

}\end{remark}

\subsection{Double-sided reflection} The objective of this section is to show that Thm.\ \ref{MT} carries over to finite-buffer systems. We use an argumentation very similar to the one used for the case of single-sided reflection. 

With $K>0$ denoting the system's buffer capacity, one can write \cite[Section 1]{LM}  the workload under double-sided reflection as
\[Q_t = Q_0 + X_t + R_t^- -R_t^+,\]
with the `regulators' $(R_t^-)_{t\geqslant 0}$ and $(R_t^+)_{t\geqslant 0}$ non-decreasing processes satisfying 
\[\int_0^\infty X_t\,{\rm d}R_t^- = \int_0^\infty (K-X_t)\,{\rm d}R_t^+=0.\]
In self-evident notation we thus have, for $s\leqslant t$, with $R_{s,t}^-$ representing the local time in 0 in $[s,t]$, and $R_{s,t}^+$ the local time in $K$ in $[s,t]$,
\[Q_t =Q_s+X_{s,t} +R_{s,t}^- -R_{s,t}^+.\]
Therefore, for $s\leqslant t\leqslant u$,
\[Q_u-Q_t = X_{t,u} + R_{t,u}^- - R_{t,u}^+.\]
Like in the system with single-sided reflection (presented in the previous subsection), $Q_t$ is non-decreasing in $Q_s$, and in addition $Q_u-Q_t$ is non-increasing in $Q_s$. The former statement follows by comparing two instances of the workload process, one being at $q_-\in[0,K)$ at time $s$, and the other at $q_+\in (q_-,K]$; then the claim can be concluded by observing that, when using the same driving L\'evy process after time $s$, the trajectories of both workload processes do not cross. The latter statement is effectively a consequence of $R_{t,u}^-$ (recalling that it represents the local time in 0 in $[t,u]$) being non-increasing in $Q_s$ and $R_{t,u}^+$ (i.e.,  the local time in $K$ in $[t,u]$) being non-decreasing in $Q_s$;
formally, this property can be established by using the following lemma.

\begin{lemma} Let $0\leqslant t\leqslant u$. Then $Q_u-Q_t$ is non-increasing in $Q_0.$
\end{lemma}

{\it Proof}: We first define, for a given initial level $Q_0\in[0,K]$ and a path $X$ of the L\'evy process, a reflection  from below at level 0. To this end, we introduce, as in (\ref{UR}),  the functional $\Gamma^-[\,\cdot\,]$, given by
\[\Gamma^-[Y]_t:= Y_t +\max\left\{0,-\inf_{0\leqslant u\leqslant t}Y_u\right\};\]
The single-sided reflection (from below, at level $0$) we obtain by imposing $\Gamma^-[\,\cdot\,]$ on $Y$, where $Y_s=Q_0+X_s$, in line with (\ref{UR}).

Now introduce the functional $\Gamma^{+}[\,\cdot\,]$
by
\[\Gamma^+[Y]_t :=Y_t-\sup_{0\leqslant s\leqslant t}\left(\max\{Y_s-K,0\}, \inf_{s\leqslant u\leqslant t} Y_u\right).\]
Then, according to \cite{KRU}, for our model with double reflection,
\[Q_t = (\Gamma^+ \circ \Gamma^-)[Q_0+X]_t= \Gamma^+[ \Gamma^-[Q_0+X]]_t.\]
In other words, the process $(Q_t)_{t\geqslant 0}$ is obtained by imposing $\Gamma^+\circ \Gamma^-$ on the process $Q_0+(X_t)_{t\geqslant 0}.$

Now consider paths $Y^{(1)}$ and $Y^{(2)}$ such that $Y^{(1)}_t\leqslant Y^{(2)}_t$ for all $t\geqslant 0.$ Then, considering the functional  $\Gamma^-[\,\cdot\,]$,
\begin{align*}
\Gamma^-[Y^{(2)}]_t - Y^{(2)}_t=&\:\max\left\{0,-\inf_{0\leqslant u\leqslant t}Y^{(2)}_u\right\}\\
\leqslant&\:\max\left\{0,-\inf_{0\leqslant u\leqslant t}Y^{(1)}_u\right\}= \Gamma^-[Y^{(1)}]_t - Y^{(1)}_t,
\end{align*}
implying that
\begin{equation}\label{bound1}\Gamma^-[Y^{(2)}]_t -\Gamma^-[Y^{(1)}]_t \leqslant Y^{(2)}_t-Y^{(1)}_t.\end{equation}
A similar reasoning applies for the functional $\Gamma^+[\,\cdot\,]$: for trajectories $Y^{(1)}$ and $Y^{(2)}$ such that $Y^{(1)}_t\leqslant Y^{(2)}_t$ for all $t\geqslant 0$,
\begin{align*}
\Gamma^+[Y^{(2)}]_t - Y^{(2)}_t=&-\sup_{0\leqslant s\leqslant t}\left(\max\{Y^{(2)}_s-K,0\}, \inf_{s\leqslant u\leqslant t} Y^{(2)}_u\right)\\
\leqslant&-\sup_{0\leqslant s\leqslant t}\left(\max\{Y^{(1)}_s-K,0\}, \inf_{s\leqslant u\leqslant t} Y^{(1)}_u\right)=\Gamma^+[Y^{(1)}]_t - Y^{(1)}_t,
\end{align*}
such that
\begin{equation}\label{bound2}\Gamma^+[Y^{(2)}]_t -\Gamma^+[Y^{(1)}]_t \leqslant Y^{(2)}_t-Y^{(1)}_t.\end{equation}
The next step is to combine the bounds (\ref{bound1}) and (\ref{bound2}).  Our goal is to show that $Q_t-Q_0$ is non-increasing in $Q_0$ (conditional on the path of the driving L\'evy process, say $X$).  To this end, consider again two instances of the workload process, one being at $Q_0:=q_-\in[0,K]$ at time $s$, and the other at $Q_0:=q_+\in [q_-,K]$; we thus have that the workloads (at time $t\geqslant 0$) corresponding with both instances are given by
\[(\Gamma^+ \circ \Gamma^-)[q_-+X]_t,\:\:\:\mbox{and}\:\:\:(\Gamma^+ \circ \Gamma^-)[q_++X]_t.\]
First recall that $\Gamma^-[q_-+X]_t \leqslant \Gamma^-[q_++X]_t$. 
Hence applying (\ref{bound2}) yields, for any   $t\geqslant 0$,
\begin{align*}\Gamma^+[\Gamma^-[q_++X]]_t-\Gamma^+[\Gamma^-[q_-+X]]_t\leqslant&\: \Gamma^-[q_++X]_t-\Gamma^-[q_-+X]_t.\end{align*}
Now using   (\ref{bound1}), the quantity in the right-hand side of the previous display is no larger than $q_+-q_-.$
Combining the above, we conclude
\[\Gamma^+[\Gamma^-[q_++X]]_t-q_+\leqslant \Gamma^+[\Gamma^-[q_-+X]]_t-q_-;\]
recalling that $q_+\geqslant q_-$, 
this is equivalent to $Q_t-Q_0$ being non-increasing in $Q_0$ (for a given path $X$, that is). By shifting time, this entails that, for $0\leqslant t\leqslant u$, $Q_u-Q_t$ being non-increasing in $Q_t$. Combining this with the fact that $Q_t$ is non-decreasing in $Q_0$, we have proven the claim. 
 \hfill$\Box$

\vspace{1mm}

Note that the above lemma immediately implies that (for $s\leqslant t\leqslant u$) $Q_u-Q_t$ is non-increasing in $Q_s$. This means that 
from this point on,
the argumentation we have developed for the infinite-buffer case can be applied in the finite-buffer context as well. We conclude that Thm.\ \ref{MT} extends to the finite-buffer case.

\begin{theorem}\label{MT2} In the system with the L\'evy process $(X_t)_{t\in{\mathbb R}}$ being reflected from below at 0 and from above at $K>0$, 
$r(\cdot)$ is non-negative, non-increasing, and convex.
\end{theorem}

In the setting above there is reflection from below at $0$ and  from above at $K>0$, but evidently this strip can be shifted. In this way we can cover the setting of reflection from below at $K_1$ and from above at $K_2$ (with $K_2>K_2$).

\section{Markov additive input: model and preliminaries}
In the following sections, our objective is to analyze the workload correlation function of a queue fed by spectrally one-sided Markov additive input. In this section, we first provide a formal definition of {\sc map}s, and introduce the notation needed. Then we present a number of useful preliminary results. 

\subsection{Model and notation}
As mentioned, we study the workload of a queue fed by a {\sc map}. 
A {\sc map} is a bivariate Markovian process $(X_t,J_t)$,
defined as follows.
\begin{itemize}
\item[$\circ$]
Let $(J_t)_t$ be an irreducible
continuous-time Markov chain with finite state space
$E=\{1,\ldots,d\}$. Define by $\DQ:=(q_{ij})_{i,j=1}^d$
the $(d\times d)$ transition rate matrix
of $(J_t)_t$. Let
${\boldsymbol \pi}$ be the (unique) stationary distribution; recall that
${\bs \pi}^{\rm T}\DQ={\bs 0}^{\rm T}$. In addition, $q_{ij}\geqslant 0$ if $i\not=j$ and
 $\DQ{\bs 1} = {\bs 0}$ (i.e., the row sums are 0). Define $\hat q _i:=-q_{ii}\geqslant 0.$
\item[$\circ$]
For each state $i$ that $J_t$ can attain, let
the process
$(X^{(i)}_t)_t$ be a L\'evy process. We either assume that either all of these are {\sc sp} or that
all of them are {\sc sn}. In the former case they have Laplace exponents
\[\varphi_i(\alpha):=\log\left({\mathbb E}\, \exp(-\alpha X^{(i)}_1)\right),\] 
whereas in the latter case they have cumulant generating functions
\[\Phi_i(\beta):=\log\left({\mathbb E}\, \exp(\beta X^{(i)}_1)\right),\] 
for $i=1,\ldots,d$
\item[$\circ$]
Letting
$T_n$ and $T_{n+1}$ be two successive transition epochs of $J_t$,
and given that $J_t$ jumps from state $i$ to state $j$ at $t=T_n$,
we define the additive process $X_t$ for $t$ in the time interval
$[T_n,T_{n+1})$ through
\[
X_t:=X_{T_n-}+U_{ij}^n+[X^{(j)}_t-X^{(j)}_{T_n}],
\]
where the $(U_{ij}^n)_n$ constitute a sequence of i.i.d.\  random
variables (each of which is distributed as a generic random variable $U_{ij}$),
describing the jumps at transition epochs. In the {\sc sp} case we work with the
{Laplace-Stieltjes} transform (with $i\in\{1,\ldots,d\}$)
\[
b_{ij}(\alpha) = {\mathbb E}\, {\rm e}^{-\alpha U_{ij}},\] where  $U_{ij}\geqslant 0$ almost surely.
Likewise, in the {\sc sn} case
we work with the moment generating functions
\[
B_{ij}(\beta) = {\mathbb E}\, {\rm e}^{\beta U_{ij}},\] 
where $U_{ij}\leqslant 0$ almost surely.
\end{itemize}
Define, in the {\sc sp} case, 
\[\sigma_{ij}(\alpha,t) :={\mathbb E}\big({\rm e}^{-\alpha X_t}1_{\{J_t=j\}}\,|\,J_0=i\big),\:\:\:
s_{ij}(\alpha):= 1_{\{i=j\}}\big( \varphi_i(\alpha)+\hat q_i\big)+ 1_{\{i\not=j\}}q_{ij} b_{ij}(\alpha),\]
and in the {\sc sn} case,
\[\Sigma_{ij}(\beta,t) :={\mathbb E}\big({\rm e}^{\beta X_t}1_{\{J_t=j\}}\,|\,J_0=i\big),
\:\:\:
S_{ij}(\beta):= 1_{\{i=j\}} \big(\Phi_i(\beta)+\hat q_i\big)+1_{\{i\not=j\}}q_{ij} B_{ij}(\beta).\]
Then, as pointed out in e.g.\ \cite[Prop. 11.3]{DM}, the matrices with entries $\sigma_{ij}(\alpha,t)$ and 
$\Sigma_{ij}(\beta,t)$ can be written in terms of matrix exponentials. 
More concretely, 
\[\sigma(\alpha,t)={\rm e}^{s(\alpha)\,t},\:\:\:
 \Sigma(\beta,t)={\rm e}^{S(\beta)\,t},\]
 for a matrix  $s(\alpha)$ with entries $s_{ij}(\alpha)$, and a matrix   $S(\beta)$ with entries $S_{ij}(\beta)$.
 
Throughout we write  ${\mathbb P}_i(A) := {\mathbb P}(A\,|\, J_0=i)$ and ${\mathbb E}_i(Y) := {\mathbb E}(Y\,|\, J_0=i)$ for events $A$ and random variables $Y$. Also, $\langle{\bs a},{\bs b}\rangle := {\bs a}^{\rm T}{\bs b}$ for ${\bs a}, {\bs b}\in{\mathbb R}^d.$
 
 The queueing processes are assumed to be stable, meaning that we impose, for the {\sc sp} and {\sc sn} case, respectively,
 \[-\langle {\bs\pi}, {\bs \varphi}'(0)\rangle +\sum_{i=1}^d \sum_{j=1}^d\pi_i q_{ij} \cdot{\mathbb E}\,U_{ij}<0,\:\:\:\:\:\:\langle {\bs\pi}, {\bs \Phi}'(0)\rangle +\sum_{i=1}^d\sum_{j=1}^d \pi_i  q_{ij} \cdot{\mathbb E}\,U_{ij}<0.\]

\subsection{Objectives and preliminary results} We now define the workload covariance function through its transform with respect to time, as follows. 
Assume the queueing process is in stationarity at time $0$.
The goal is to study, for $\vt>0$,
\[\gamma(\vt) := {\mathbb C}{\rm ov}(Q_0,Q_T),\]
where $T\equiv T(\vt)$ is exponentially distributed with mean $\vt^{-1}$, independently of the driving {\sc map}.  
Withe $r(t)$ the correlation between $Q_0$ and $Q_t$ and $\hat r(\cdot)$ the Laplace transform of $r(\cdot)$, we thus have the evident relation
\begin{equation}\label{TR}\hat r(\vt) =\frac{\gamma(\vt)}{\vt}\cdot \frac{1}{{\mathbb V}{\rm ar}\, Q_0}.\end{equation}
In the next two sections we focus on deriving an expression for $\gamma(\vt)$ (for the {\sc sp} and {\sc sn} cases, respectively), whereas the examples in Section \ref{maie} provide insight in the shape of $r(\cdot)$ and its asymptotics for $t$ large.

If the underlying queueing system would have been Markovian, this expression could have been rewritten as
\[{\mathbb E}(Q_0 Q_T) - ({\mathbb E}\,Q_0)^2 =
\int_{0}^{\infty} x \,m(x) \,p(x)\,{\rm d}x- ({\mathbb E}\,Q_0)^2 ,\]
with  $p(x)\,{\rm d}x:=\,{\mathbb P}(Q_0\in {\rm d}x)$ and $m(x) := {\mathbb E}(Q_T\,|\,Q_0=x)$.
In our case, however $(Q_t)_{t\geqslant 0}$ (alone) is {\it not} Markovian, but $(Q_t, J_t)_{t\geqslant 0}$ jointly are. Define $p_i(x) = {\mathbb P}(Q_0\in {\rm d}x,J_0=i)$, and $m_i(x):=  {\mathbb E}(Q_T\,|\,Q_0=x,J_0=i).$ In evident vector notation, we obtain that in our setting with {\sc map} input
\begin{equation}\label{DEF}\gamma(\vt) = \int_{0}^{\infty} x \,\langle {\bs m}(x) ,{\bs p}(x)\rangle \,{\rm d}x- ({\mathbb E}\,Q_0)^2 .\end{equation}
The following result will be used several times. Define $\tau(x):=\inf\{t\geqslant 0: X_t\leqslant -x\}$ for $x\geqslant 0$.

\begin{lemma} \label{L1} For all $x\geqslant 0$,
\[m_i(x)=  x +\int_x^\infty {\mathbb E}_i\,{\rm e}^{-\vt\,\tau(y)}\,{\rm d}y + {\mathbb E}_iX_T.\]
\end{lemma}
{\it Proof}:
The quantity ${\bs m}(x)$ can be determined by investigating the effect of slightly perturbing the initial workload level $x$. More precisely, we compare $m_i(x+\delta)$ and $m_i(x)$ for $\delta$ small. Observe that if before $T$ the queues have drained, the queueing processes have coupled, and as a consequence the workloads are the same; if, on the contrary, the busy period has not ended before $T$, the workloads still differ by $\delta.$
This entails that, as $\delta\downarrow 0$, for an $f\in(0,1)$,
\begin{align*}m_i(x+\delta)-m_i(x) = \:&\delta \,{\mathbb P}_i(\tau(x) > T)+f\delta \, \,{\mathbb P}_i(\tau(x) < T<\tau(x+\delta))\\
= \:&\delta \,{\mathbb P}_i(\tau(x) > T)+f\delta \, \,{\mathbb P}_i(\underline X_T\in (-x-\delta,-x))\\
= \:&\delta \,{\mathbb P}_i(\tau(x) > T)+o(\delta)
,\end{align*}
with $\underline X_t$ denoting the minimum of $X_s$ over $s\in [0,t]$.
Dividing both sides of the previous display by $\delta$ and letting $\delta\downarrow 0$, we thus obtain the differential equation
\[m_i'(x) = {\mathbb P}_i(\tau(x) > T).\]
By conditioning on 
the value of $\tau(x)$, 
\[
 {\mathbb P}_i(\tau(x) > T) = 1- {\mathbb E}_i\,{\rm e}^{-\vt\,\tau(x)}.\]
We thus find that, for  constants $C$ and $K$,
\[m_i(x)= x -\int_0^x {\mathbb E}_i\,{\rm e}^{-\vt\,\tau(y)}\,{\rm d}y + C =  x +\int_x^\infty {\mathbb E}_i\,{\rm e}^{-\vt\,\tau(y)}\,{\rm d}y + K.\]
The constant $K$ can be determined as follows. As $x\to\infty$, we have that $m_i(x)- x \downarrow {\mathbb E}_iX_T$ (the underlying reason being that as the initial level $x$ increases to $\infty$, the probability of reflection before $T$ vanishes). It thus follows that $K= {\mathbb E}_iX_T.$  $\hfill\Box$

\vb

Later in this paper we use a result for the transform of ${\mathbb E}_i\,{\rm e}^{-\vt\,\tau(y)}$ (with respect to $y$, that is). To this end, recall that we have observed above that
\begin{align*}{\mathbb E}_i\,{\rm e}^{-\vt\,\tau(y)} = {\mathbb P}_i(\tau(y)< T) 
.\end{align*}Recall that ${\mathbb P}_i(\tau(y)< t) = {\mathbb P}_i(\underline X_t < -y).$
Upon combining these two facts,
\[\int_0^\infty{\rm e}^{-\eta y}\,{\mathbb E}_i\,{\rm e}^{-\vt\,\tau(y)}\,{\rm d}y = 
\int_0^\infty {\rm e}^{-\eta y}\,{\mathbb P}_i(\underline X_T < -y)\,{\rm d}y.\]
Now using integration by parts, we end up with the following result. We remark that its validity is not restricted to spectrally one-sided L\'evy processes. 
\begin{lemma} \label{L2} For $\eta\geqslant 0$,
\[\int_0^\infty{\rm e}^{-\eta y}\,{\mathbb E}_i\,{\rm e}^{-\vt\,\tau(y)}\,{\rm d}y =\frac{1}{\eta}\left(1-
{\mathbb E}_i\,{\rm e}^{\eta \underline{X}_T}\right).\]
\end{lemma}

\section{Markov additive input: spectrally positive case}
The main objective of this section is to identify the transform $\gamma(\cdot)$ for the case of {\sc sp-map} input.  The main idea is to use representation (\ref{DEF}). To do so, we first collect a number of useful expressions.

In the first place, in this {\sc sp-map} case the Laplace transform of $\kappa(\alpha)$ of $Q_0$ is given by
\[\kappa(\alpha) :={\mathbb E}\,{\rm e}^{-\alpha Q_0}= \alpha\, {\bs \ell}\, (s(\alpha))^{-1} \,{\bs 1};\]
see for instance \cite[Eqn. (11.1)]{DM} and \cite[Corollary 3.1]{DiM}. 
The  row vector ${\bs \ell}$ can be evaluated using  \cite[Lemmas 3.1 and 3.2]{DiM}.
This result can be seen as the extension of the classical Pollaczek-Khinchine formula \cite[Thm.~3.2]{DM} for spectrally positive L\'evy input \cite{ZO} to the spectrally positive {\sc map} case. 
It follows that ${\mathbb E}\,Q_0 = -\kappa'(0).$ In addition,
${\mathbb V}{\rm ar}\,Q_0= \kappa''(0) - (\kappa'(0))^2.$

A somewhat more precise form of the above result for $\kappa(\alpha)$ concerns the distribution of $Q_0$ jointly with the background state $J_0$:
\[\kappa_i(\alpha) :={\mathbb E}\,({\rm e}^{-\alpha Q_0}\,1_{\{J_0=i\}})= \alpha\left( {\bs \ell} \,(s(\alpha))^{-1} \right)_i.\]

Bearing in mind Eqn.\ (\ref{DEF}) and Lemma \ref{L1}, the next question is: how to evaluate the quantities $\xi_i(y,\vt):= {\mathbb E}_i\,{\rm e}^{-\vt\,\tau(y)}$ and
$\nu_i\equiv\nu_i(\vt):={\mathbb E}_iX_T$; if we could evaluate these, we have an expression for $\gamma(\vt).$ 

\begin{itemize}
\item[$\circ$]
We can determine $\nu_i\equiv\nu_i(\vt)$ as follows.
As we have seen before
\[{\mathbb E}_i\,{\rm e}^{-\alpha X_t} = \big({\rm e}^{s(\alpha)\,t}\,{\bs 1}\,\big)_i,\]
so that
\[\nu_i={\mathbb E}_i X_T = \int_0^\infty \vt \,{\rm e}^{-\vt t}\left(-\lim_{\alpha\downarrow 0} 
\frac{{\rm d}}{{\rm d}\alpha} \big({\rm e}^{s(\alpha)\,t}\,{\bs 1}\,\big)_i\right){\rm d}t.\]
There is, however, an easier way to compute $\nu_i$. 
Observe that, using the memoryless property, 
\[\nu_i = - \frac{\varphi_i'(0)}{\hat q_i +\vt} +
\sum_{j\not=i} \frac{q_{ij}}{\hat q_i+\vt} \left({\mathbb E}\,U_{i} + \nu_j\right);\]
this (diagonally dominant) linear system has a unique solution, and can be solved by standard methods.
\item[$\circ$]
We now point out how to evaluate $ {\mathbb E}_i\,{\rm e}^{-\vt\,\tau(y)}$. The following standard reasoning can be used.
In the first place, observe that $(\tau(y), J_{\tau(y)})_{y\geqslant 0}$ is a {\sc map}; $\tau_x$ is increasing (and therefore {\sc sp}), and (due to the drift     condition) $\tau_x\uparrow \infty$ as $x\to\infty.$
Now define the matrix $A(\vt)$ whose $(i,j)$-th component is given by
\[A_{ij}(\vt):= {\mathbb E}( {\rm e}^{-\vt\,\tau(0+) }1_{\{J_{\tau(0+) =j\}}}|\,J_0=i). \]
Because $(\tau(x), J_{\tau(x)})_{x\geqslant 0}$ is an {\sc sp-map}, there is matrix-valued function $\Psi(\cdot)$ such that
\[{\mathbb E}({\rm e}^{-\vt\,\tau(y) }1_{\{J_{\tau(y) =j\}}}|\,J_{0+}=i) = ({\rm e}^{-\Psi(\vt)\,y})_{i,j};\]
for characterizations of $\Psi(\cdot)$, see e.g.\ \cite{DM,KP}. We thus conclude that there are matrices $A(\vt)$ and $\Psi(\vt)$ such that
\[ \xi_i(y,\vt):= {\mathbb E}_i\,{\rm e}^{-\vt\,\tau(y)} = \big(A(\vt)\,{\rm e}^{-\Psi(\vt) \,y}\,{\bs 1}\,\big)_i.\]
We refer for more background to e.g.\ \cite{KP}, and,  for the case of phase-type jumps, \cite{BR}.
Depending on the linear-algebraic properties of the matrix $\Psi(\vt)$, this expression can be further evaluated. For instance if all eigenvalues of $\Psi(\vt)$ are distinct (say $\psi_1(\vt)$ up to $\psi_d(\vt)$), we can write
\begin{equation}
\label{expa} \xi_i(y,\vt)= \sum_{j=1}^d c_{ij}(\vt) \, {\rm e}^{-\psi_j(\vt)\,y}\end{equation}
for  appropriate $c_{ij}(\vt)$.
For the moment we assume that this expansion applies; if it does not apply (when eigenvalues have multiplicities larger than 1), then straightforward adaptations can be made, as explained in Remark \ref{R2}.
\end{itemize}

We thus obtain that
\[m_i(x) = x +\int_x^\infty  \sum_{j=1}^d c_{ij}(\vt) \, {\rm e}^{-\psi_j(\vt)\,y}\, {\rm d}y +\nu_i(\vt)=
x+\sum_{j=1}^d \frac{c_{ij}(\vt)}{\psi_j(\vt)} {\rm e}^{-\psi_j(\vt)\,x} +\nu_i(\vt).\]
We are now in a position to evaluate $\gamma(\vt).$ We have
\begin{align*}\gamma(\vt) = &\,\int_{0}^{\infty} x \cdot\sum_{i=1}^d \left(x+\sum_{j=1}^d \frac{c_{ij}(\vt)}{\psi_j(\vt)} {\rm e}^{-\psi_j(\vt)\,x} +\nu_i(\vt)\right) p_i(x){\rm d}x - ({\mathbb E}\,Q_0)^2 \\
=&\:{\mathbb V}{\rm ar}\, Q_0 -\sum_{i=1}^d \sum_{j=1}^d \frac{c_{ij}(\vt)}{\psi_j(\vt)}\kappa_i'(\psi_j(\vt))+\sum_{i=1}^d \nu_i(\vt) \pi_i \,{\mathbb E}_i Q_0.
\end{align*}
Here $\pi_i \,{\mathbb E}_i Q_0 = -\kappa_i'(0),$ and, as was mentioned earlier, ${\mathbb V}{\rm ar}\,Q_0= \kappa''(0) - (\kappa'(0))^2.$ We have thus found the following representation for $\gamma(\vt)$. 

\begin{theorem}\label{SP}
If the expansion $(\ref{expa})$ applies, then, for any $\vt\geqslant 0$,
\[\gamma(\vt) =\kappa''(0) - (\kappa'(0))^2 -\sum_{i=1}^d \sum_{j=1}^d \frac{c_{ij}(\vt)}{\psi_j(\vt)}\kappa'_i(\psi_j(\vt))-\sum_{i=1}^d \nu_i(\vt) \kappa_i'(0).\]
\end{theorem}

\begin{remark}{\em The above theorem can be considered as the {\sc map} counterpart of the corresponding result for a queue driven by a spectrally positive L\'evy process \cite[Thm. 2.1]{EM}. Observe that, so as to translate the results into one another, the expression for $\gamma(\vt)$ has to be divided by $\vt$ as well as by ${\mathbb V}{\rm ar}\,Q_0$, cf.\ (\ref{TR}).}\hfill $\diamond$
\end{remark}

\begin{remark}\label{R2}{\em 
Condition (\ref{expa}) is imposed for convenience. If there are mixed exponential-polynomial terms in the expression for $\xi_i(y,\vt)$ (due to eigenvalues with multiplicities higher than 1), one could proceed as follows. Consider  for instance the case that the eigenvalue $\psi_j(\vt)$ having multiplicity two leads to a term proportional to $y\, {\rm e}^{-\psi_j(\vt)\,y}.$ Due to
the identity
\[\int_x^\infty y\,{\rm e}^{-\psi y}\,{\rm d}y = \frac{x}{\psi}\,{\rm e}^{-\psi x} +\frac{1}{\psi^2}\,{\rm e}^{-\psi x},\]
we find that for this $j$ (and all $i\in\{1,\ldots,d\}$) the $-\kappa'_i(\psi_j(\vt))/\psi_j(\vt)$ appearing in Thm.\  \ref{SP} has to be replaced by
\[
\frac{\kappa''(\psi_j(\vt))}{\psi_j(\vt)} - \frac{\kappa'_i(\psi_j(\vt))}{(\psi_j(\vt))^2}.
\]
Analogously, if a term $y^k\, {\rm e}^{-\psi_j(\vt)\,y}$ appears for some $k\in\{0,\ldots,d-1\}$, then the $(k+1)$-st derivative of $\kappa(\cdot)$ will end up in the expression for $\gamma(\vt).$
}\hfill $\diamond$
\end{remark}

\section{Markov additive input: spectrally negative case}
In this section we provide an expression for $\gamma(\vt)$ for the case of {\sc sn-map} input. As it turns out, now the quantities $\xi_i(y,\vt)= {\mathbb E}_i\,{\rm e}^{-\vt\,\tau(y)}$ can not be given (unlike in the {\sc sp-map} case), but their transforms with respect to $y$ {\it can} be evaluated. Combining this with the fact that $Q_0$ is of phase type facilitates the computation of $\gamma(\vt).$ Below we detail this reasoning.

\begin{lemma}\label{LPT}
In the {\sc sn-map}  case, $Q_0$ (jointly with the state $J_0$ being $i$) has a phase-type distribution. In particular, we can write 
$Q_0$ (jointly with $J_0=i$) is given by the initial distribution ${\boldsymbol t}_i$ and transition rate matrix $E$.
\end{lemma}

{\it Proof}:
This property follows directly using (for instance) the results from \cite{DM}. The reasoning consists of the following steps. (i)~From \cite[Prop. 4.2]{DM} it follows how to translate (by using time-reversal) the distribution of $Q_0$ (jointly with $J_0=i$) into the distribution of the all-time maximum of the non-reflected process (given that the background state at time $0$ is $i$). (ii)~Then we use the material of \cite[Section 3]{DM} on the all-time minimum for an {\sc sp-map}; clearly, one obtains results for the all-time maximum of an {\sc sn-map} by looking at the negative of the process. (iii) This means that we can make use of \cite[Thm. 3.2]{DM} that provides the transform of the all-time minimum $\underline Y$ for an {\sc sp-map} $(Y_t)_{t\geqslant 0}$ (or, alternatively, we could use
\cite[Section 6]{IV}). Note that this transform also covers the epoch the minimum is attained, but (as we are interested in $\underline Y$ only) one has to set the corresponding parameter in the transform equal to $0$. (iv) We thus obtain that, with $\underline J$ the state of the background process at the epoch the maximum is attained, the matrix of which the $(i,j)$-th element is given by
\[{\mathbb E} ({\rm e}^{\beta\underline{Y}} 1{\{\underline J=j\}}|\,J_0= i)\]
can be written as 
\[-F (\beta I-E)^{-1} E,\]
for a matrix $F$ with non-negative elements and row-sums equal to 1, and a defective transition rate matrix $E$ (and $I$ defining the $d$-dimensional identity matrix). (v) Observe that, by e.g.  \cite[Prop.\ III.4.1(iii)]{AS}, this defines a phase-type distribution.
\hfill$\Box$

Due to Lemma \ref{LPT}, by \cite[Prop.\ III.4.1(i)]{AS}, 
\[{\mathbb P}(Q_0\in{\rm d}x, J_0=i) = -\big( {\boldsymbol t}_i^{\rm T} \,{\rm e}^{Ex} E{\boldsymbol 1}\big)\,{\rm d}x.\]
For the moment we assume that, for vectors ${\bs \zeta}$ and ${\bs \eta}$, the above expression allows the expansion
\begin{equation}\label{expand}
p_i(x)\,{\rm d}x={\mathbb P}(Q_0\in{\rm d}x, J_0=i) = \sum_{j=1}^d \zeta_{ij} {\rm e}^{-\eta_j\, x}\,{\rm d}x,\end{equation}
a sufficient condition being that the eigenvalues of $E$ be distinct; in Remark \ref{multi} we comment on adaptations when (\ref{expand}) does not apply. In a way the above assumption can be considered as the {\sc sn}-counterpart of the assumption (\ref{expa}) that we imposed in the {\sc sp-}case. It follows that
\[\kappa_i(\alpha):={\mathbb E}({\rm e}^{-\alpha Q_0}1\{J_0=i\}) =\sum_{j=1}^d\zeta_{ij} \frac{1}{\eta_j+\alpha},\]
so that
\[\kappa_i'(0) = -{\mathbb E}_iQ_0 = \sum_{j=1}^d\frac{\zeta_{ij}}{\eta_j},\:\:\:\:\kappa''(0)-(\kappa'(0))^2={\mathbb V}{\rm ar}\, Q_0 = 
2\sum_{j=1}^d\frac{\zeta_{ij}}{\eta_j^2}-\left(\sum_{j=1}^d\frac{\zeta_{ij}}{\eta_j}\right)^2.\]

Recall from Lemma \ref{L2} that the transform of ${\mathbb E}_i{\rm e}^{-\vt\,\tau(y)}$ (with respect to $y$, that is) can be expressed in terms of the transform of $\underline X_T$.
In addition, we remark that an expression for this transform of $\underline X_T$, which we denote
by 
\[\Omega_i(\eta):= {\mathbb E}_i {\rm e}^{\eta\underline X_T},\]
is given in e.g.\ \cite[Thm.\ 1.(ii), Eqn.\ (17)]{KlP} and \cite[Prop. 6.1]{IV}. Now, similarly to what we have done  in the {\sc sp}-case, $\gamma(\vt)$ can be obtained appealing to (\ref{DEF}) and Lemma \ref{L1}. To this end, we have that
\begin{align*}
\int_0^\infty x\,\langle {\bs m}(x),{\bs p}(x)\rangle\, {\rm d}x&=\,\sum_{i=1}^d \int_0^\infty x\left(x+\int_x^\infty{\mathbb E}_i{\rm e}^{-\vt\,\tau(y)}{\rm d}y +{\mathbb E}_i X_T\right) \sum_{j=1}^d \zeta_{ij} {\rm e}^{-\eta_j\, x}{\rm d}x\\
&=\,{\mathbb E}(Q_0^2)+ \sum_{i=1}^d  \sum_{j=1}^d\int_0^\infty \int_x^\infty x\,{\mathbb E}_i{\rm e}^{-\vt\,\tau(y)} \zeta_{ij} {\rm e}^{-\eta_j\, x}{\rm d}y\,{\rm d}x + \sum_{i=1}^d \nu_i(\vt) {\mathbb E}_i Q_0,
\end{align*}
with $\nu_i(\vt)$ as defined (and evaluated) in the {\sc sp}-case. With (recalling Lemma \ref{L2})
\[\bar\Omega_i(\eta,\vt):=\frac{1}{\eta}\big(1-\Omega_i(\eta)\big)=\int_0^\infty {\rm e}^{-\eta y}  \,{\mathbb E}_i{\rm e}^{-\vt\,\tau(y)} {\rm d}y  ,\]
we have (by swapping the order of the integrals)
\begin{align*}
\int_0^\infty \int_x^\infty x\,{\mathbb E}_i{\rm e}^{-\vt\,\tau(y)} \zeta_{ij} {\rm e}^{-\eta_j\, x}{\rm d}y\,{\rm d}x\,&=\int_0^\infty {\mathbb E}_i{\rm e}^{-\vt\,\tau(y)} \zeta_{ij}\int_0^y x{\rm e}^{-\eta_j\, x}{\rm d}x\,{\rm d}y\\
\,&=\frac{\zeta_{ij}}{\eta_j^2}\int_0^\infty {\mathbb E}_i{\rm e}^{-\vt\,\tau(y)} \left(1-{\rm e}^{-\eta_j y}-\eta_jy\,{\rm e}^{-\eta_j y}\right){\rm d}y\\
\,&=\frac{\zeta_{ij}}{\eta_j^2}\left(\bar\Omega_i(0,\vt)-\bar\Omega_i(\eta_j,\vt)+\eta_j\bar\Omega'_i(\eta_j,\vt)\right),
\end{align*}
where $\bar\Omega'_i(\eta,\vt)$ is understood as the derivative of $\bar\Omega_i(\eta,\vt)$ with respect to $\eta$. Subtracting $({\mathbb E}Q_0)^2$, we thus arrive at the following result.

\begin{theorem} If the expansion $(\ref{expand})$ applies, then, for any $\vt\geqslant 0$,
\[\gamma(\vt) =\kappa''(0) - (\kappa'(0))^2 -\sum_{i=1}^d \sum_{j=1}^d\frac{\zeta_{ij}}{\eta_j^2}\left(\bar\Omega_i(0,\vt)-\bar\Omega_i(\eta_j,\vt)+\eta_j\bar\Omega'_i(\eta_j,\vt)\right)-\sum_{i=1}^d \nu_i(\vt) \kappa_i'(0)
.\]
\end{theorem}

\begin{remark}\label{multi}{\em 
Regarding condition $(\ref{expand})$, we remark that
one may have situations in which mixed exponential-polynomial terms show up in $p_i(x)$ (due to non-simple eigenvalues). These can be dealt with as in Remark \ref{R2}.} \hfill $\diamond$
\end{remark}

\section{Markov additive input: examples}\label{maie}
In this section we treat a number of illustrative examples. In the first example we present explicit computations for one of  the most frequently used {\sc map}s, viz.\ Markov-modulated fluid. The second example is such that the structural properties, that have been derived for the L\'evy case ($r(\cdot)$ being non-negative, non-increasing, convex), still apply, as well as the tail behavior agreeing with the tail of the busy-period distribution. The third example then shows that in general $r(\cdot)$ is {\it not} non-negative, non-increasing, convex; in addition the rate of decay of $r(t)$ for $t$ large is shown to potentially differ from that of the busy-period tail, depending on properties of the modulating process $(J_t)_{t\in{\mathbb R}}$. 
\subsection{Markov-modulated fluid with a general number of states}
We define the Markov-modulated fluid model as follows. In any of the states of the background process, the process behaves as a (state-specific) drift. 
Let $\mu_i$ be the constant drift when the background state is~$i$, where $i\in\{1,\ldots,d\}$; assume that the $\mu_i$\,s do not equal $0$. Let ${\mathscr E}_-$ the subset of $\{1,\ldots,d\}$ such that $\mu_i<0$ and 
let ${\mathscr E}_+$ the subset of $\{1,\ldots,d\}$ such that $\mu_i>0$.
We start by considering the case of a general dimension $d\in\{2,3,\ldots\}.$ We first explain how ${\bs \nu}(\vt)$ and ${\bs \xi}(y,\vt)$ can be found.
\begin{itemize}
\item[$\circ$]
It is straightforward to determine the mean rate $\nu_i(\vt).$ It takes some elementary algebra to verify that ${\bs \nu}(\vt) = -(\DQ-\vt I)^{-1}{\bs \mu}.$

\item[$\circ$]
We set up a system of linear differential equations for the $\xi_i(y,\vt)$. Using standard arguments, as $\delta\downarrow 0$,
\begin{align*}
\xi_i(y,\vt) =&\: \sum_{j\not=i} q_{ij}\delta \,\xi_j(y,\vt) + \left(1-\hat q_i \delta\right) e^{-\vt\delta} \xi_i(y+\mu_i\delta,\vt)+o(\delta)\\
=&\: \sum_{j\not=i} q_{ij}\delta \,\xi_j(y,\vt) + \left(1-\hat q_i \delta-\vt\delta\right) \xi_i(y+\mu_i\delta,\vt) +o(\delta).
\end{align*}
It thus follows that, where $\xi'(y,\vt)$ corresponds to differentiation with respect to $y$,
\[-\mu_i\,\xi'_i(y,\vt) = \sum_{j=1}^d q_{ij}\,\xi_j(y,\vt) -\vt\,\xi_i(y,\vt).\]
In self evident notation, this is, in vector/matrix form,
\[{\bs \xi}'(y,\vt) =-\Psi(\vt)\,{\bs \xi}(y,\vt),\]
with $\Psi(\vt):={\mathscr M}^{-1}(\DQ- \vt I)$ and ${\mathscr M}:={\rm diag}\{{\bs \mu}\}.$ By solving the system of linear differential equations, it thus follows that ${\bs \xi}(y,\vt) ={\rm e}^{-\Psi(\vt)\,y}{\bs \xi}(0+,\vt),$
where
\[ {\xi}_i(0+,\vt) :=\lim_{y\downarrow 0} \xi_i(y,\vt).\]
Observe that ${\xi}_i(0+,\vt)=1$ for $i\in{\mathscr E}_-$, whereas ${\xi}_i(0+,\vt)<1$ for $i\in{\mathscr E}_+$; in the latter case ${\xi}_i(0+,\vt)$ is the transform of the queue's busy period given the background process is  in state $i$ at the start of this busy period.
\end{itemize}
\subsection{Markov-modulated fluid with two states} In this case the steady-state distribution follows from the results of \cite{KK}. Write $q_1:=\hat q_1 = q_{12}$, $q_2:=\hat q_2 =q_{21}$, and $\bar q:=q_1+q_2$. The steady-state distribution of the background process is ${\bs \pi} = (q_2/\bar q,q_1/\bar q).$ Following \cite{KK}, we assume that $\mu_1>0$ and $\mu_2<0$. The average drift is
\[\bar\mu := \pi_1 \mu_1+\pi_2\mu_2 =\frac{\mu_1 q_2}{\bar q} +\frac{\mu_2 q_1}{\bar q},\]
which we throughout assume to be negative. 
We define $\lambda:= q_2/\mu_1 +q_1/\mu_2>0.$
As argued in \cite{KK},
\[p_1(x) = \pi_1 \lambda \,{\rm e}^{-\lambda x}, \:\:\:p_2(x) := \pi_2 \frac{q_2}{q_2-\lambda\mu_2}
\lambda \,{\rm e}^{-\lambda x}.\]
This entails that
\begin{align*}\kappa_1(\alpha) =&\, \int_0^\infty {\rm e}^{-\alpha x}  \pi_1 \lambda \,{\rm e}^{-\lambda x}\,{\rm d}x
= \pi_1 \frac{ \lambda}{\lambda +\alpha},\\
\kappa_2(\alpha)=&\,\int_0^\infty {\rm e}^{-\alpha x} \pi_2 \frac{q_2}{q_2-\lambda\mu_2}
\lambda \,{\rm e}^{-\lambda x}\,{\rm d}x=\pi_2 \frac{ \lambda}{\lambda +\alpha}\,\frac{q_2}{q_2-\lambda\mu_2}
.
\end{align*}
Hence,
\[-\kappa'_1(0)=\pi_1\,{\mathbb E}_1Q_0=\frac{\pi_1}{\lambda},\:\:\:
-\kappa'_2(0)=\pi_2\,{\mathbb E}_2Q_0=\frac{\pi_2}{\lambda}\frac{q_2}{q_2-\lambda\mu_2}
\]
and
\[\kappa_1''(0) = \frac{2\pi_1}{\lambda^2},\:\:\: \frac{2\pi_2}{\lambda^2}\frac{q_2}{q_2-\lambda\mu_2},\:\:\:\:{\mathbb E}(Q_0^2)=\frac{2}{\lambda^2}\left(\pi_1+\pi_2\frac{q_2}{q_2-\lambda\mu_2}\right).\]
We now evaluate ${\bs \nu}(\vt)$:
\[\left(\begin{array}{c}
\nu_1(\vt)\\\nu_2(\vt)\end{array}\right) = \left(\begin{array}{cc}\vt+q_1&-q_1\\-q_2&\vt+q_2
\end{array}\right)^{-1}\left(\begin{array}{c}
\mu_1\\ \mu_2\end{array}\right) =\frac{1}{\vt^2+\vt \bar q}\left(\begin{array}{c}(\vt+q_2)\mu_1+q_1\mu_2\\q_2\mu_1+(\vt+q_1)\mu_2
\end{array}\right).\]
We now concentrate on determining $\xi_1(y,\vt)$ and $\xi_2(y,\vt).$ In this special case there is a direct approach; we start by analyzing $\xi_2(y,\vt).$ Let $N(x)$ be Poisson distributed with mean $x>0$. Let $\tau$ denote a busy period (starting at workload level 0, with the background process just having turned to state 1), and $\tau_1,\tau_2,\ldots$ i.i.d.\ copies of $\tau$; also, $\pi(\vt) :={\mathbb E}\,{\rm e}^{-\vt\tau}$.
 Then, using basic self-similarity properties, given that $J_0=2$,
\[\tau(y)\stackrel{\rm d}{=}- \frac{y}{\mu_2} +\sum_{i=1}^{N(-q_2y/\mu_2)} \tau_i.\]
In other words, with $\varrho(\vt):=-(\vt +q_2(1- \pi(\vt)))/\mu_2$,
\[\xi_2(y,\vt) ={\rm e}^{-\vt y/\mu_2} \sum_{k=0}^\infty {\rm e}^{-q_2 y/\mu_2} \frac{(q_2 y/\mu_2)^k}{k!} \,\big({\mathbb E}\,{\rm e}^{-\vt\tau}\big)^k={\rm e}^{-\varrho(\vt)\,y}.\]
It is not hard to see that
\[\xi_1(y,\vt) = \pi(\vt) \xi_2(y,\vt)=\pi(\vt)\,{\rm e}^{-\varrho(\vt)\,y}.\]
So we are left with determining $\pi(\vt)$. This we do by relating $\tau$ to the busy period $\tau^\circ$ in an M/M/1 queue with arrival rate $\alpha$, service rate $\beta$, and depletion rate $R$. With $\pi^\circ(\vt)\equiv\pi^\circ(\vt\,|\,\alpha,\beta,R):={\mathbb E}\,{\rm e}^{-\vt\tau^\circ}$, 
we have
\[\pi^\circ(\vt) = \int_0^\infty {\beta}\,{\rm e}^{-\beta x}\,{\rm e}^{-\vartheta x/R}
\sum_{k=0}^{\infty}{\rm e}^{-\alpha x/R}\frac{(\alpha x/R)^k}{k!} \big(\pi^\circ(\vt)\big)^k\,{\rm d}x=\frac{\beta R}{\beta R +\vt+\alpha(1-\pi^\circ(\vt))}.\]
Now $\pi^\circ(\vt)$ can be found be solving the above second-order equation and picking the correct root. A straightforward geometric argument then yields that
\[\pi(\vt) =\pi^\circ\left(\vt\left.\left(1-\frac{\mu_2}{\mu_1}\right)\,\right|\,q_2,\frac{q_1}{\mu_1},-\mu_2\right).\]
It then follows that
\begin{align*}
m_1(x)=&\:x+\frac{\pi(\vt)}{\varrho(\vt)}\,{\rm e}^{-\varrho(\vt)\,x} +\frac{\mu_1\vartheta+\bar \mu \bar q}{\vt^2+\vt \bar q},\:\:\:\:
m_2(x)=\: x +\frac{1}{\varrho(\vt)}\,{\rm e}^{-\varrho(\vt)\,x} +\frac{\mu_2\vartheta+\bar \mu \bar q}{\vt^2+\vt \bar q}.
\end{align*}
Now $\gamma(\vt)$ can be explicitly evaluated. Define $c_1(\vt):=\pi(\vt)$, $c_2(\vt):=1$, and
\[a_1 = \pi_1 \lambda, \:\:\:a_2 := \pi_2 \frac{q_2}{q_2-\lambda\mu_2}
\lambda.\]
It thus follows that
\begin{align*}\gamma(\vt)=&\: \sum_{i=1}^2 \int_0^\infty x \left(x+{c_i(\vt)} 
{\rm e}^{-\varrho(\vt)\,x}+ \frac{\mu_i\vartheta+\bar \mu \bar q}{\vt^2+\vt \bar q}\right) 
a_i \,{\rm e}^{-\lambda x}{\rm d}x-\left(\frac{a_1+a_2}{\lambda^2}\right)^2\\
=&\:
 \sum_{i=1}^2 a_i\left(\frac{2}{\lambda^3}
 -\frac{c_i(\vt)}{(\varrho(\vt)+\lambda)^2}+
  \frac{\mu_i\vartheta+\bar \mu \bar q}{\vt^2+\vt \bar q}\,\frac{1}{\lambda^2}
 \right)-\left(\frac{a_1+a_2}{\lambda^2}\right)^2
.\end{align*}
For ease we renormalize time and space such that $\mu_1=q_2=1$; we simply put $\mu:=-\mu_2>0$ and $q:=q_1.$ Hence,
\[\pi(\vt) =\pi^\circ(\vt(1+\mu)\,|\,1,q,\mu),\]
which equals
\[\pi(\vt) =\frac{1}{2}\left(q\mu +\vt(1+\mu)+1 - \sqrt{(q\mu +\vt(1+\mu)+1)^2-4q\mu}\right).\]
The asymptotic behavior of $c(t):={\mathbb C}{\rm ov}(Q_0,Q_t)$ for $t$ large can be found by inspecting 
the singularities of $\gamma(\vt).$ More specifically, as described in full detail in \cite[Section 3]{AW}, the rightmost singularity of $\gamma(\vt)$ in the left half-plane determines the tail behavior of $c(t)$. In our case there is a pole at $\vt_{\rm p}=-\bar q=-q-1$. In addition, there is a branching point where $\pi(\vt)$ has a branching point, which is at
\[\vt_{\rm b} =\frac{2\sqrt{q\mu}-q\mu-1}{1+\mu} = -\frac{(\sqrt{q\mu}-1)^2}{1+\mu}.\]
Noting that $2\sqrt{q\mu} \geqslant -\mu-q$,
\[\vt_{\rm b} =\frac{2\sqrt{q\mu}-q\mu-1}{1+\mu}  \geqslant\frac{-\mu-q-q\mu-1}{1+\mu} = -q-1 =\vt_{\rm p}.\]
We conclude that in this example the tail behavior of $c(t)$ is determined by $\vt_{\rm b}.$ More specifically,
\[\lim_{t\to\infty} \frac{1}{t} \log c(t) = \vt_{\rm b}.\]

\subsection{Example illustrating tail behavior of correlation}
In this example the background process has  a cyclic structure: $q_{12} = q_{23} =\cdots=q_{d-1,d}=q_{d1} =q$,
$\mu_1>0$, and $\mu_i=-\infty$ for $i\in\{2,\ldots,d\}.$ 
To make the notation more compact, we again renormalize time and space by putting $\mu_1=q=1.$
Strictly speaking, this example is not a {\sc map} (due to $\mu_i=-\infty$ for $i\in\{2,\ldots,d\}$), but it can be approximated arbitrarily closely by a {\sc map}. By this example we show that the tail behavior of $c(\cdot)$ does not necessarily match that of the tail of the busy period distribution, depending on specific features of the matrix $\DQ.$

As it turns out, in this setting we can compute ${\bs p}(x)$ and ${\bs m}(x)$ explicitly. It is easily seen that busy periods are exponentially distributed with mean $1$ and that during such busy periods the background state is~1. When the background process is in $i\in\{2,\ldots,d\}$, the workload level is $0$. Hence $p_1(x) = {\rm e}^{-x}/d$, whereas ${\mathbb P}(Q_0=0,J_0=i) = 1/d$ for $i\in\{2,\ldots,d\}$.
Also,
\[m_1(x) = \frac{\vt}{1+\vt}\,\left(x+\frac{1}{1+\vt}\right)+
\frac{1}{1+\vt}\,{\mathbb P}(E_{d-1} <T) \,m_1(0),
\]
where $E_{d-1}$ denotes an Erlang random variable with $d-1$ phases, each phase being exponentially distributed with parameter 1. As
\[{\mathbb P}(E_{d-1} <T)  = \left(\frac{1}{1+\vt}\right)^{d-1},\]
we conclude that 
\[m_1(0)=\frac{\vt(1+\vt)^{d-2}}{(1+\vt)^d-1},\:\:\:m_1(x) = \frac{\vt}{1+\vt}\,\left(x+\frac{1}{1+\vt}\right)+
\left(\frac{1}{1+\vt}\right)^d\,
\frac{\vt(1+\vt)^{d-2}}{(1+\vt)^d-1}
.\]
It thus follows that
\[\gamma(\vt) = \int_0^\infty \frac{x}{d}\,{\rm e}^{-x} \left(\frac{\vt}{1+\vt}\,\left(x+\frac{1}{1+\vt}\right)+
\left(\frac{1}{1+\vt}\right)^d\,
\frac{\vt(1+\vt)^{d-2}}{(1+\vt)^d-1}
\right) {\rm d}x-\frac{1}{d^2},\]
which simplifies to
\begin{align}\nonumber\gamma(\vt)=&\: \frac{2}{d}\frac{\vt}{1+\vt} +\frac{1}{d} 
\frac{\vt}{(1+\vt)^2}+\frac{1}{d}
\frac{1}{(1+\vt)^2}\,
\frac{\vt}{(1+\vt)^d-1}-\frac{1}{d^2}\\
=&\: \frac{2}{d}\frac{\vt}{1+\vt} +\frac{1}{d} \left(\frac{\vt(1+\vt)^{d-2}-d^{-1}((1+\vt)^d-1)}{(1+\vt)^d-1}\right)
\label{gamm}.\end{align}
This expression has poles at $\vt=-1$ and at $\vt=-1+{\rm e}^{2\pi{\rm i}\,k/d}$, for $k \in\{1,2,\ldots, d-1\}$; observe that,  in the expression between the brackets in (\ref{gamm}), the zero in the denominator for $\vt=0$ is compensated by a zero in the corresponding numerator, entailing that $\vt=0$ is not a pole.  
We now show that increasing $d$ drastically affects the nature of the asymptotics of $c(t)$.

\begin{itemize}
\item[$\circ$]
We start with $d=2.$ Then
\begin{align*}\gamma(\vt) =&\:\frac{\vt}{1+\vt} +\frac{1}{2} 
\frac{\vt}{(1+\vt)^2}  +\frac{1}{2}\frac{1}{(1+\vt)^2}\frac{1}{2+\vt}
-\frac{1}{4}\\
=&\:\frac{3}{4}-\frac{1}{2}\frac{1}{1+\vt}-\frac{1}{2}\frac{1}{(1+\vt)^2} +
\frac{1}{2}
\frac{1}{(1+\vt)^2}\frac{1}{2+\vt}.\end{align*}
This we rewrite to 
\begin{align*}\int_0^\infty {\rm e}^{-\vt t}c(t){\rm d}t =&\: \frac{3}{4\vt}-\frac{1}{2\vt}\left(\frac{1}{1+\vt}\right)-\frac{1}{2\vt}\left(\frac{1}{1+\vt}\right)^2
+\frac{1}{4\vt}\left(\frac{1}{1+\vt}\right)^2\left(\frac{2}{2+\vt}\right)\\
=&\:\frac{1}{\vt}\left(\frac{3}{4}-\frac{1}{1+\vt}+\frac{1}{2}\frac{1}{2+\vt}\right)=
\frac{1}{1+\vt}-\frac{1}{4}\frac{1}{2+\vt}.
\end{align*}
Explicit inversion yields
\[c(t) = {\rm e}^{-t} -\tfrac{1}{4} \,{\rm e}^{-2t}.\]
In this case $\gamma(\vt)$ has poles at $\vt=-1$ and $\vt=-2.$ As mentioned, the right-most pole dominates the asymptotic behavior. We have $c(t)\, {\rm e}^{t} \to 1$ as $t\to\infty.$
\item[$\circ$] We continue with $d=3$. Now, besides the usual pole at $\vt=-1$, there are two more poles at $\vt=
{\rm e}^{2\pi{\rm i}/3} = (-3+{\rm i}\,\sqrt{3})/2$ and ${\rm e}^{4\pi{\rm i}/3} = (-3-{\rm i}\,\sqrt{3})/2$, both of them having real part $-3/2$. By comparing these  poles, it follows that the pole with the largest real part remains $\vt=-1$.
We have
\[\int_0^\infty {\rm e}^{-\vt t}c(t){\rm d}t =
\frac{2}{3}\frac{1}{1+\vt}
-\frac{1}{9}\frac{\vt}{\vt^2+3\vt+3},\]
eventually leading to
\[c(t) = \frac{2}{3}{\rm e}^{-t} +\frac{\sqrt{3}}{9} {\rm e}^{-3t/2}\sin(\tfrac{1}{2}\sqrt{3}\,t)
-\frac{1}{9} {\rm e}^{-3t/2}\cos(\tfrac{1}{2}\sqrt{3}\,t),\]
and thus $c(t)\, {\rm e}^{t} \to 2/3$ as $t\to\infty.$ Conclude that we see essentially the same behavior as for $d=2.$
\item[$\circ$] For $d=4$ we have
\[\int_0^\infty {\rm e}^{-\vt t}c(t){\rm d}t = \frac{1}{2}\frac{1}{1+\vt}
-\frac{1}{16}\frac{\vt^2-2}{\vt^3+4\vt^2+6\vt+4}.\]
In this case the poles are, besides $\vt=-1$, also $\vt=-2,$ $-1+{\rm i}$, and $-1-{\rm i}.$ After tedious but straightforward calculus,
\[c(t) =\frac{1}{2} {\rm e}^{-t} -\frac{1}{16}{\rm e}^{-2t} +\frac{1}{8}{\rm e}^{-t}\sin t,\]
so that $\lim_{t\to\infty} c(t)\,{\rm e}^{-t}$ does not exist. More specifically,
\[\liminf_{t\to\infty} c(t)\,{\rm e}^{-t} = \frac{3}{8},\:\:\:\:\limsup_{t\to\infty} c(t)\,{\rm e}^{-t} = \frac{5}{8}.\]
Conclude that in this case there are three poles with the same `magnitude' (in terms of the size of their real parts),
which consequently all appear in the asymptotic behavior. 
\item[$\circ$] We now consider $d\in \{5,6,\ldots\}$. The crucial observation is that, for $d$ larger than $5$,  the poles  $-1 +{\rm e}^{2\pi{\rm i}\,/d}$
and  $-1 +{\rm e}^{2\pi{\rm i}\,/d}$ are dominant, as these have a real part that is strictly larger than $-1$; this real part is $\phi_d:= -1+ \cos(2\pi/d).$
As in the case of $d=4$, $\lim_{t\to\infty} c(t)\,{\rm e}^{-t}$ does not exist because of oscillating terms. Logarithmic asymptotics can be derived, though: the decay rate $\bar c$ equals
\[\bar c:=\lim_{t\to\infty}\frac{1}{t}\log c(t) = \phi_d.\]
\end{itemize}
So for $d$ up to 4 the decay rate corresponds with the pole at $-1$, whereas for larger values of $d$ the pole at $\cos(2\pi/d)$ takes over: we have found that $\bar c =-1 +\max\{0,\cos(2\pi/d)\}.$

\begin{remark}\label{coupling}{\em The fact that the decay rate of $c(\cdot)$ does not necessarily coincide with that of the tail of the busy period distribution can be understood as follows.  In \cite{GM} a coupling-based approach was developed for the case of L\'evy input. The argumentation relied on writing $c(t)$ as ${\mathbb E}(Q_0(Q_t -Q\s_t))$, with the process $(Q\s_t)_{t\geqslant 0}$ being constructed as follows: sample $Q_0\s$ independently from $Q_0$, but use the same driving Levy process $(X_t)_{t\geqslant 0}$ as for $(Q_t)_{t\geqslant 0}$.  It is now easily seen that a condition for $c(t)$ to be non-zero is that at least one of the two busy periods has not ended at time $t$; more precisely, we should have that $\underline X_t>-Q_0$ or $\underline X_s>-Q_0\s$. This explains why in the L\'evy case the tail behavior of $c(\cdot)$ is intimately connected with the tail distribution of the busy period. These ideas have been made precise in \cite{GM}. 

In the case of {\sc map} input however, the states of the background processes corresponding to the sample  $Q_0$ and $Q_0\s$ do not necessarily coincide. As a consequence, we cannot use the same driving L\'evy process $(X_t)_{t\geqslant 0}$ after time 0, and hence the coupling idea does not apply. To remedy this complication, an idea could be to let the L\'evy processes $(X_t)_{t\geqslant 0}$ and $(X\s_t)_{t\geqslant 0}$ run independently till the background processes reach the same state, and apply the above coupling from that point on. It may, however, take relatively long before the background states meet.
It is precisely this insight that explains why the structure of the background process $(J_t)_{t\in{\mathbb R}}$ (represented by its transition rate matrix $\DQ$) plays a role. This is nicely illustrated in the above example: when the dimension $d$ grows, the coupling takes longer, entailing that from a certain $d$ on (in the example $d=4$) the tail of $c(\cdot)$ starts to diverge from the tail of the busy-period distribution. } \hfill $\diamond$
\end{remark}

\section*{Appendix}
{\it Proof of Lemma \ref{LEM}}: We believe this result has appeared, in various forms, in the literature. Several approaches can be followed to prove the claim; for completeness we include a compact proof here. The starting point is the relation
\begin{equation}
\label{AB}{\mathbb E}(AB) = \int_0^\infty\int_0^\infty {\mathbb P}(A>a, B>b)\,{\rm d}a\,{\rm d}b,\end{equation}
for any non-negative $A$ and $B$ for which the above objects are well-defined. 

First suppose that $f(\cdot)$ is non-negative; then, by (\ref{AB}),
\begin{align*}{\mathbb E}(A\,f(A)) &\: =  \int_0^\infty\int_0^\infty {\mathbb P}(A>a, f(A)>b)\,{\rm d}a\,{\rm d}b,\\
{\mathbb E}(B\,f(A)) &\: =  \int_0^\infty\int_0^\infty {\mathbb P}(B>a, f(A)>b)\,{\rm d}a\,{\rm d}b.\end{align*}
Then observe that
\begin{align*}
{\mathbb P}(B>a, f(A)>b) \leqslant&\: \min\{{\mathbb P}(B>a), {\mathbb P}(f(A)>b)\}\\
=&\:\min\{{\mathbb P}(A>a), {\mathbb P}(f(A)>b)\}={\mathbb P}(A>a, f(A)>b),
\end{align*}
using the fact that $f(\cdot)$ is non-decreasing in the (non-strict) inequality, and
\[{\mathbb P}(A>a, f(A)>b)={\mathbb P}(A>\max\{a, f^{-1}(b)\})=\min\{{\mathbb P}(A>a), {\mathbb P}(A>f^{-1}(b))\} \]
in the last equality. We thus establish the claim under the proviso that $f(\cdot)$ is non-decreasing and non-negative.

We now lift the condition that $f(\cdot)$ be non-negative. To this end, suppose that $f(0)<0$. Then, realizing that $a\mapsto f(a)-f(0)$ is non-negative, by applying the previous result,
\begin{align*}
{\mathbb E}(A\,f(A))&\:= {\mathbb E}\big(A\,(f(A)-f(0))\big) -f(0)\,{\mathbb E}(A) \\
&\:\geqslant {\mathbb E}\big(B\,(f(A)-f(0))\big) -f(0)\,{\mathbb E}(A) = {\mathbb E}(B\,f(A)).
\end{align*}
The result for $g(\cdot)$ follows analogously. \hfill$\Box$

\bibliographystyle{plain}

\end{document}